# AUTOMATIC DETECTION OF THE COMMON AND NON- COMMON FREQUENCIES IN CONGRUENT DISCRETE SPECTRA – A THEORETICAL APPROACH

Cezar DOCA, and Constantin PAUNOIU

*Pitesti, Romania*

cezar_doca@yahoo.com

ABSTRACT

*Both sampling a time-varying signal, and its spectral analysis are activities subjected to theoretically compelling, such as Shannon's theorem and the objectively limiting of the frequency's resolution. Usually, the signals' spectra are processed and interpreted by a scientist who, presumably, has sufficient prior information about the monitored signals to conclude on the significant frequencies, for example. On the other hand, processing and interpretation of signals' spectra can be routine tasks that must be automated using suitable software (PC application). In the above context, the paper presents the theoretic bases of an intuitive and practical approach of the (automatic) detection of the common and non-common frequencies in two or more congruent spectra.*

## 1. Introduction

Our paper is neither about time-varying signal's sampling and/or his spectral analysis, these specialized subjects being fully analyzed in a lot of very good theoretical works as [1] – [10], nor about some possible algebraic structures and/or other involved mathematical aspects applicable on a set of congruent spectra.

We only theoretically develop (and propose) a practical method dedicated to help the experimenter in a quickly (and automatically) identification of the common and non-common frequencies existing in two ore more congruent discrete spectra. A method which uses, to conclude on the significant frequencies, an intuitive discernment criterion based on the magnitudes of the spectral lines, namely: *the more important frequency in signal has the greater magnitude among the spectral lines.*

## 2. Four definitions

If $s(t)$ is a time-varying signal, then:

**D1.** *We imagine the signal's* ***(frequencies) discrete spectrum*** *as a collection (a set)* $S(\Delta\nu) = \{(\nu_i, A(\nu_i)) | i = 0,1,\cdots,N\}$ *of* $N+1$ *doublets* $(\nu, A(\nu))$, *where* $A(\nu)$ *is the magnitude of the spectral line having as abscissa the frequency ν; particularly, the discrete spectrum of a time-constant signal* $s_{ct}(t) = const.$ *is* $S_{ct}(\Delta\nu) = \{(0, A_{ct}(0))\}$;

**D2.** *The value* $\Delta\nu = \nu_i - \nu_{i-1} = const.$ *is the discrete spectrum's* ***resolution****, and* $\nu_i = i \cdot \Delta\nu$;

**D3.** *We name two discrete spectra* $S_1(\Delta\nu_1)$ *and* $S_2(\Delta\nu_2)$ *having the same resolution, i.e.* $\Delta\nu_1 = \Delta\nu_2 = \Delta\nu$, *as congruent*;

**D4.** *For two congruent discrete spectra* $S_1(\Delta\nu)$ *and* $S_2(\Delta\nu)$ *we name the spectral lines having the same abscissa, i.e.* $(\nu_k, A_1(\nu_k)) \in S_1(\Delta\nu)$ *and* $(\nu_k, A_2(\nu_k)) \in S_2(\Delta\nu)$, *as* ***correspondent****.*

Starting from the above unconventional definition of the discrete spectrum $S(\Delta\nu)$ we observe that, if all $A(\nu_i)\neq 0$, condition usually accomplished in practical measurements, then:

**O1.** *Noting* $A^{-1}(\nu_i)=1/A(\nu_i)$, *the set* $S^{-1}(\Delta\nu)=\{(\nu_i, A^{-1}(\nu_i)) | i=0,1,\ldots,N\}$ *is a discrete spectrum congruent with* $S(\Delta\nu)$ *and, in accordance with the intuitive discernment criterion based on the magnitudes of the spectral lines, the frequencies significant in* $S(\Delta\nu)$ *are insignificant in* $S^{-1}(\Delta\nu)$ *and vice versa*; *symbolically we write* $S^{-1}(\Delta\nu)=1/S(\Delta\nu)$.

**O2.** *If* $S_1(\Delta\nu)$, $S_2(\Delta\nu)$, ..., $S_M(\Delta\nu)$ *are congruent discrete spectra then, noting* $A_{\alpha\beta\cdots\zeta}(\nu_i)=A_\alpha(\nu_i)\cdot A_\beta(\nu_i)\cdot\ldots\cdot A_\zeta(\nu_i)$, $\alpha,\beta,\ldots\zeta\in\{1,2,\ldots M\}$, *all the possible collections* $S_{\alpha\beta\cdots\zeta}(\Delta\nu)=\{(\nu_i, A_{\alpha\beta\cdots\zeta}(\nu_i)) | i=0,1,\ldots,N\}$ *are also discrete spectra congruent with* $S_\alpha(\Delta\nu)$, $S_\beta(\Delta\nu)$, ..., $S_\zeta(\Delta\nu)$; *symbolically we write* $S_{\alpha\beta\cdots\zeta}(\Delta\nu)=S_\alpha(\Delta\nu)\cdot S_\beta(\Delta\nu)\cdot\ldots\cdot S_\zeta(\Delta\nu)$.

## 3. Two propositions

Let $s_1(t)$ and $s_2(t)$ be two time-varying signals; if $S_1(\Delta\nu)$ and $S_2(\Delta\nu)$ are their congruent discrete spectra, then the next two propositions are true:

**P1.** *The discrete spectrum* $S_p(\Delta\nu)=S_1(\Delta\nu)\cdot S_2(\Delta\nu)$ *contains as emphasized the* ***common frequencies correspondent lines****, significantly existing in both* $S_1(\Delta\nu)$ *and* $S_2(\Delta\nu)$ *spectra*; *certainly,* $S_p(\Delta\nu)$ *is congruent with* $S_1(\Delta\nu)$ *and* $S_2(\Delta\nu)$.

Demonstration: accepting the spectral line $(\nu_k, A_p(\nu_k))\in S_p(\Delta\nu)$ as a particular case of the general situation presented in **O2**, it is clear that the product value $A_p(\nu_k)=A_1(\nu_k)\cdot A_2(\nu_k)$ amplifies, most of all, when both values $A_1(\nu_k)$ and $A_2(\nu_k)$ increase, i.e. when the frequency $\nu_k$ is significant in both spectra.

**P2.** *The discrete spectrum* $S_r(\Delta\nu)=S_1(\Delta\nu)\cdot S_2^{-1}(\Delta\nu)$ *contains as emphasized the* ***non-common frequencies correspondent lines****, existing significantly in* $S_1(\Delta\nu)$ *spectrum, and insignificantly in* $S_2(\Delta\nu)$ *spectrum*; *certainly,* $S_r(\Delta\nu)$ *is congruent with* $S_1(\Delta\nu)$ *and* $S_2(\Delta\nu)$.

Demonstration: accepting the spectral line $(\nu_k, A_r(\nu_k))\in S_r(\Delta\nu)$ as a generalization of the particular case presented in **O1**, it is clear that the ratio value $A_r(\nu_k)=A_1(\nu_k)/A_2(\nu_k)$ amplifies, most of all, when the numerator value $A_1(\nu_k)$ increases, i.e. the frequency $\nu_k$ is significant in $S_1(\Delta\nu)$ spectrum, and when the denominator value $A_2(\nu_k)$ decreases, i.e. the frequency $\nu_k$ is insignificant in $S_2(\Delta\nu)$ spectrum.

## 4. Generalization for more spectra

Using the propositions **P1 – P2**, and the observations **O1 – O2**, it is easy to demonstrate that:

**P3.** *The discrete spectrum* $S_{\alpha\beta\cdots\zeta}(\Delta\nu)=S_\alpha(\Delta\nu)\cdot S_\beta(\Delta\nu)\cdot\ldots\cdot S_\zeta(\Delta\nu)$ *contains as emphasized the* ***common frequencies correspondent lines****, significantly existing in all* $S_\alpha(\Delta\nu)$, $S_\beta(\Delta\nu)$,... $S_\zeta(\Delta\nu)$ *congruent discrete spectra.*

Demonstration: according to proposition **P1**, the discrete spectrum $S_{\alpha\beta}(\Delta\nu) = S_{\alpha}(\Delta\nu) \cdot S_{\beta}(\Delta\nu)$ contains as emphasized the common frequencies correspondent lines for $S_{\alpha}(\Delta\nu)$ and $S_{\beta}(\Delta\nu)$; now, the discrete spectrum $S_{\alpha\beta\gamma}(\Delta\nu) = S_{\alpha\beta}(\Delta\nu) \cdot S_{\gamma}(\Delta\nu)$ contains as emphasized the common frequencies correspondent lines for $S_{\alpha\beta}(\Delta\nu)$ and $S_{\gamma}(\Delta\nu)$, i.e. for $S_{\alpha}(\Delta\nu)$, $S_{\beta}(\Delta\nu)$ and $S_{\gamma}(\Delta\nu)$; and so forth; finally, the discrete spectrum $S_{\alpha\beta\ldots\zeta}(\Delta\nu) = S_{\alpha\beta\ldots}(\Delta\nu) \cdot S_{\zeta}(\Delta\nu)$ contains as emphasized the common frequencies correspondent lines for $S_{\alpha\beta\ldots}(\Delta\nu)$ and $S_{\zeta}(\Delta\nu)$, i.e. the common frequencies correspondent lines for $S_{\alpha}(\Delta\nu)$, $S_{\beta}(\Delta\nu)$, ... $S_{\zeta}(\Delta\nu)$.

**P4.** *If* $S_{\alpha\beta\ldots\zeta}(\Delta\nu) = S_{\alpha}(\Delta\nu) \cdot S_{\beta}(\Delta\nu) \cdot \ldots \cdot S_{\zeta}(\Delta\nu)$ *and* $S_{ab\cdots z}(\Delta\nu) = S_{a}(\Delta\nu) \cdot S_{b}(\Delta\nu) \cdot \ldots \cdot S_{z}(\Delta\nu)$ *are congruent discrete spectra, then the discrete spectrum* $S_{\alpha\beta\cdots\zeta ab\cdots z}(\Delta\nu) = S_{\alpha\beta\cdots\zeta}(\Delta\nu) \cdot S^{-1}_{ab\cdots z}(\Delta\nu)$ *contains as empha-sized the* ***non-common frequencies correspondent lines****, existing significantly in* $S_{\alpha}(\Delta\nu)$, $S_{\beta}(\Delta\nu)$,... $S_{\zeta}(\Delta\nu)$ *congruent discrete spectra, and insignificantly in* $S_{a}(\Delta\nu)$, $S_{b}(\Delta\nu)$,... $S_{z}(\Delta\nu)$ *congruent discrete spectra.*

Demonstration: according to proposition **P2**, the discrete spectrum $S_{\alpha\beta\cdots\zeta ab\cdots z}(\Delta\nu)$ contains as emphasized the common frequencies correspondent lines for $S_{\alpha\beta\cdots\zeta}(\Delta\nu)$ and $S^{-1}_{ab\cdots z}(\Delta\nu)$; but, according to observation **O1**, the significant frequencies in $S^{-1}_{ab\cdots z}(\Delta\nu)$ are not significant in $S_{ab\cdots z}(\Delta\nu)$; i.e. each emphasized frequency of $S_{\alpha\beta\cdots\zeta ab\cdots z}(\Delta\nu)$ is a non-common frequency for $S_{\alpha\beta\cdots\zeta}(\Delta\nu)$ and $S_{ab\cdots z}(\Delta\nu)$, significant in $S_{\alpha\beta\cdots\zeta}(\Delta\nu)$, and insignificant in $S_{ab\cdots z}(\Delta\nu)$; namely, all emphasized frequencies in $S_{\alpha\beta\cdots\zeta ab\cdots z}(\Delta\nu)$ are significant in all $S_{\alpha}(\Delta\nu)$, $S_{\beta}(\Delta\nu)$,... $S_{\zeta}(\Delta\nu)$ congruent discrete spectra, and insignificant in all $S_{a}(\Delta\nu)$, $S_{b}(\Delta\nu)$,... $S_{z}(\Delta\nu)$ congruent discrete spectra.

## 5. Some demonstrative examples

Let $s_1(t)$ and $s_2(t)$ be two time-varying signals, with $S_1(\Delta\nu)$ and $S_2(\Delta\nu)$ their discrete spectra as in Figures 1 – 2, numerically simulated in Microsoft Office Excel with: *sampling time* $t_s = 2.55s$, *sampling frequency* $\nu_s = 100Hz$, *number of samples* $N_s = 256$, *number of drawn spectral lines* $N_l = 128$, and *spectral resolution* $\Delta\nu = 0.392157Hz$:

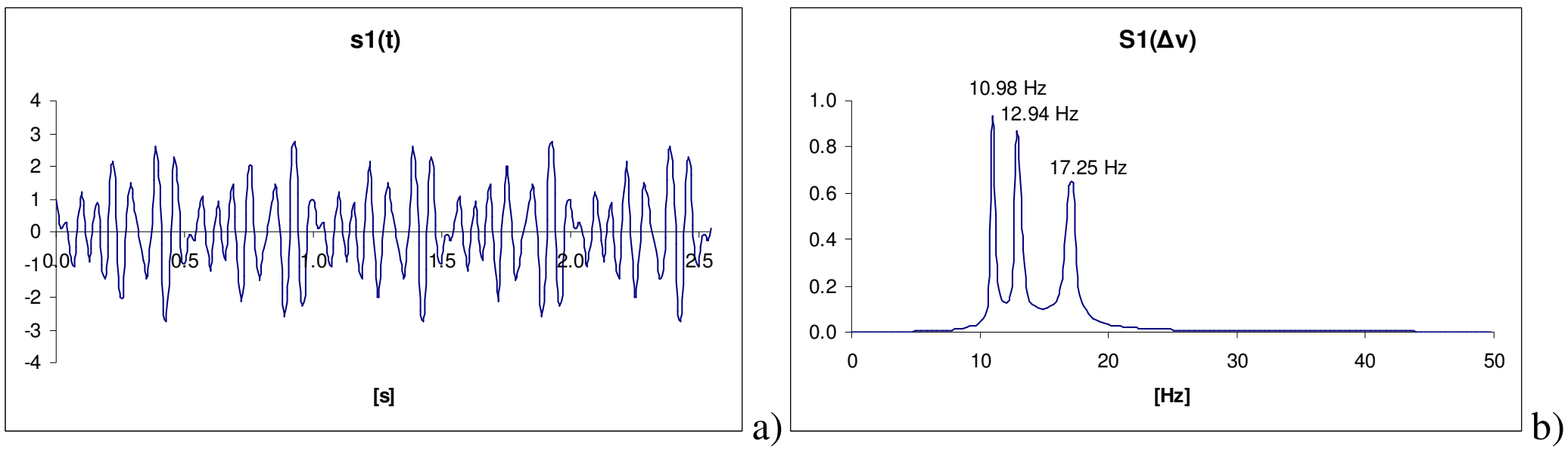

Figure 1. *Signal* $s_1(t) = \sin(2\pi \cdot 11 \cdot t) + \cos(2\pi \cdot 13 \cdot t) - \sin(2\pi \cdot 17 \cdot t)$:
*a) sampling* *b) spectrum*

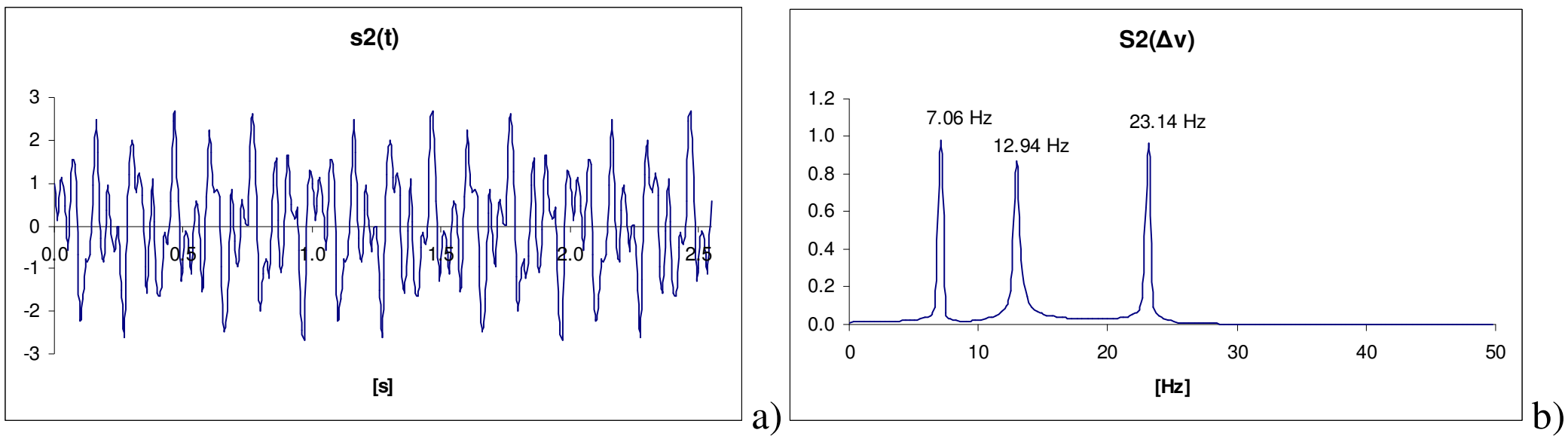

Figure 2. *Signal* $s_2(t) = \sin(2\pi \cdot 7 \cdot t) + \cos(2\pi \cdot 13 \cdot t) - \sin(2\pi \cdot 23 \cdot t)$:
*a) sampling* *b) spectrum*

Applying the proposition **P1** we obtain *the common frequencies*:

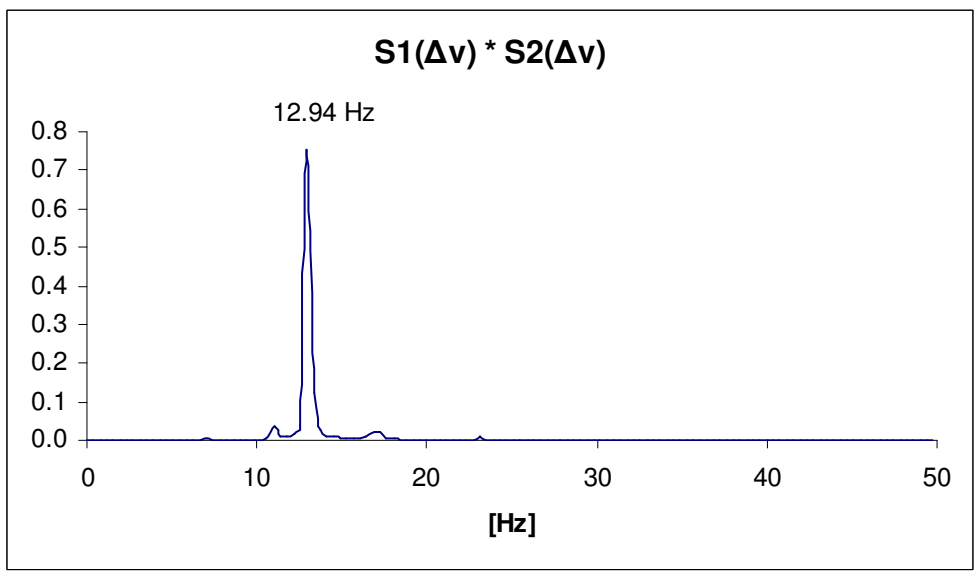

Figure 3. *Common frequencies*:
*the products* $A_1(\nu_k) \cdot A_2(\nu_k); k = 0, \ldots, 127$

Applying the proposition **P2** we obtain *the non-common frequencies*:

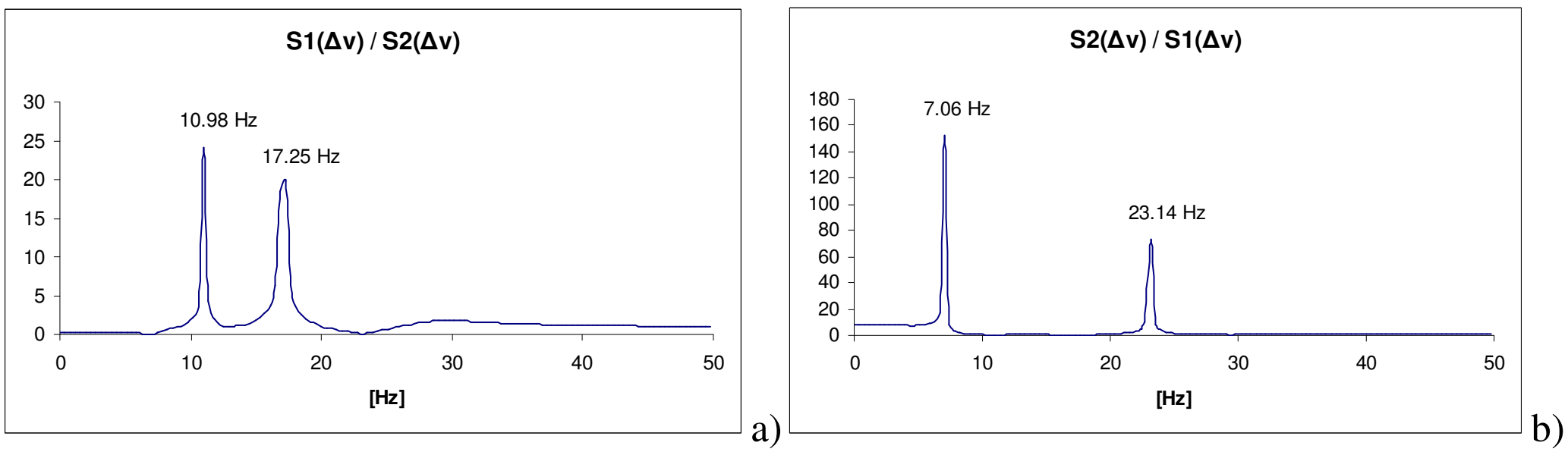

Figure 4. *Non-common frequencies*:
*a) the ratios* $A_1(\nu_k)/A_2(\nu_k); k = 0, \ldots, 127$ *b) the ratios* $A_2(\nu_k)/A_1(\nu_k); k = 0, \ldots, 127$

Let $s_3(t)$, $s_4(t)$ and $s_5(t)$ be another three time-varying signals, with $S_3(\Delta\nu)$, $S_4(\Delta\nu)$ and $S_5(\Delta\nu)$ their discrete spectra as in Figures 5 – 7 (also numerical simulations in Microsoft Office Excel with: *sampling*: $t_s = 2.55s$; $\nu_s = 100Hz$; $N_s = 256$, and *spectrum*: $N_l = 128$, $\Delta\nu = 0.392157Hz$):

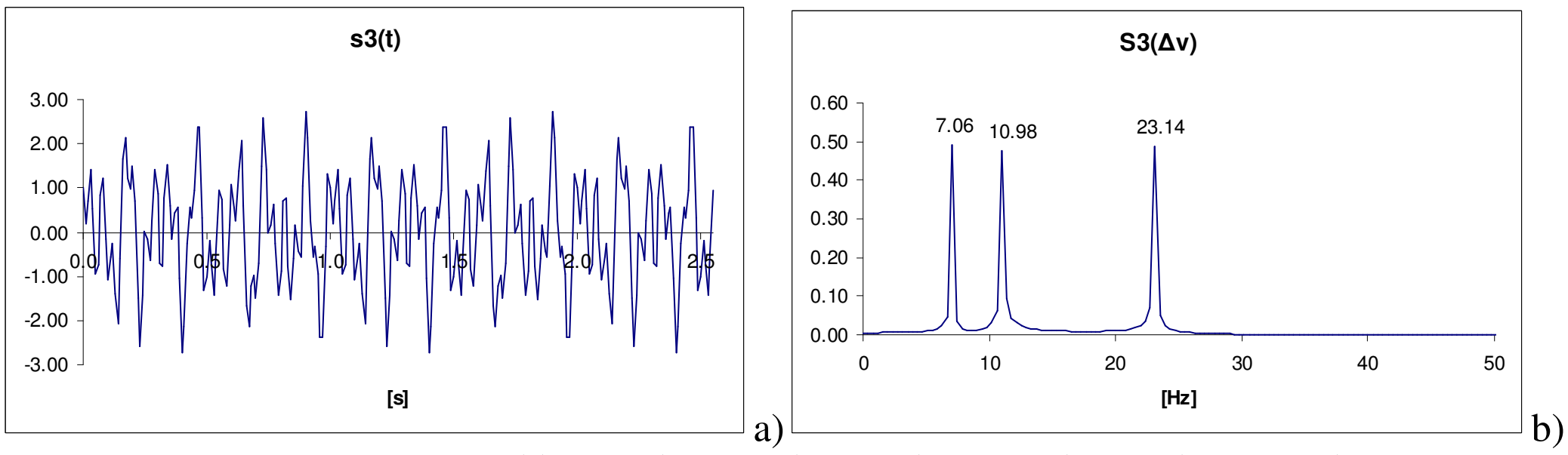

Figure 5. *Signal* $s_3(t) = \sin(2\pi \cdot 7 \cdot t) + \cos(2\pi \cdot 11 \cdot t) - \sin(2\pi \cdot 23 \cdot t)$:
*a) sampling* *b) spectrum*

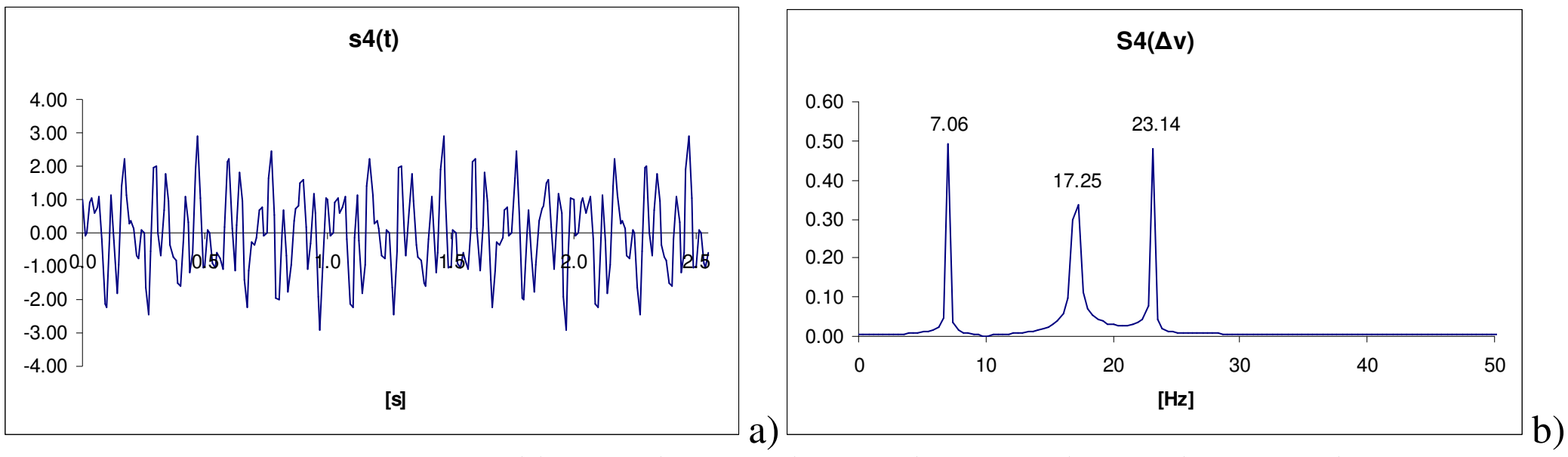

Figure 6. *Signal* $s_4(t) = \sin(2\pi \cdot 7 \cdot t) + \cos(2\pi \cdot 17 \cdot t) - \sin(2\pi \cdot 23 \cdot t)$:
*a) sampling* *b) spectrum*

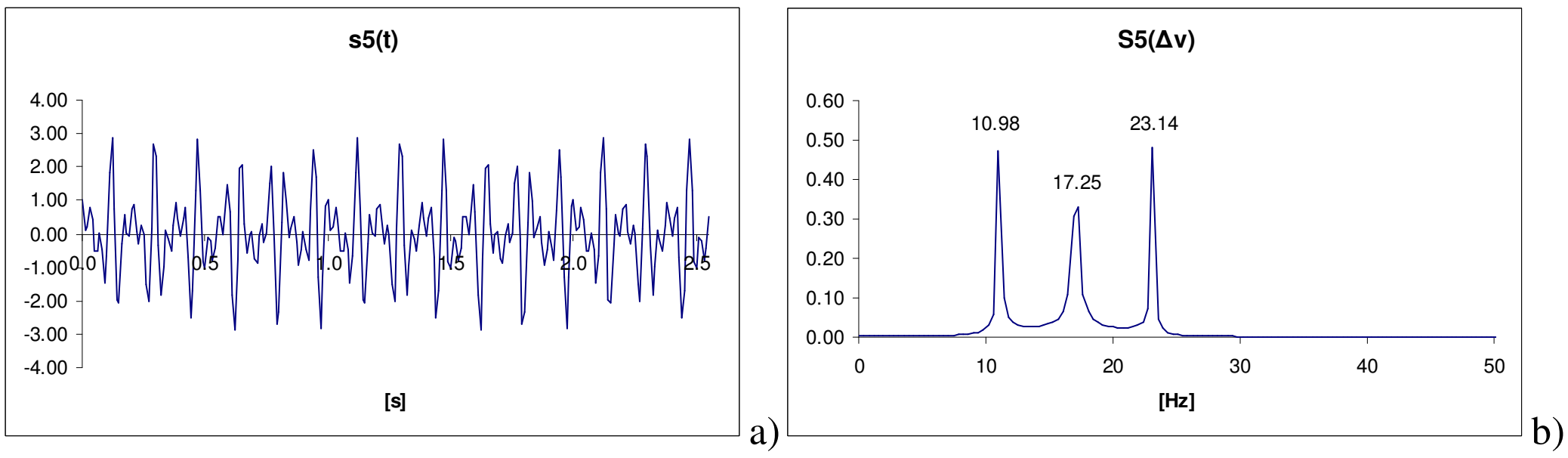

Figure 7. *Signal* $s_5(t) = \sin(2\pi \cdot 10 \cdot t) + \cos(2\pi \cdot 17 \cdot t) - \sin(2\pi \cdot 23 \cdot t)$:
*a) sampling* *b) spectrum*

Now it is easy to verify the imagined situations from the next figures:

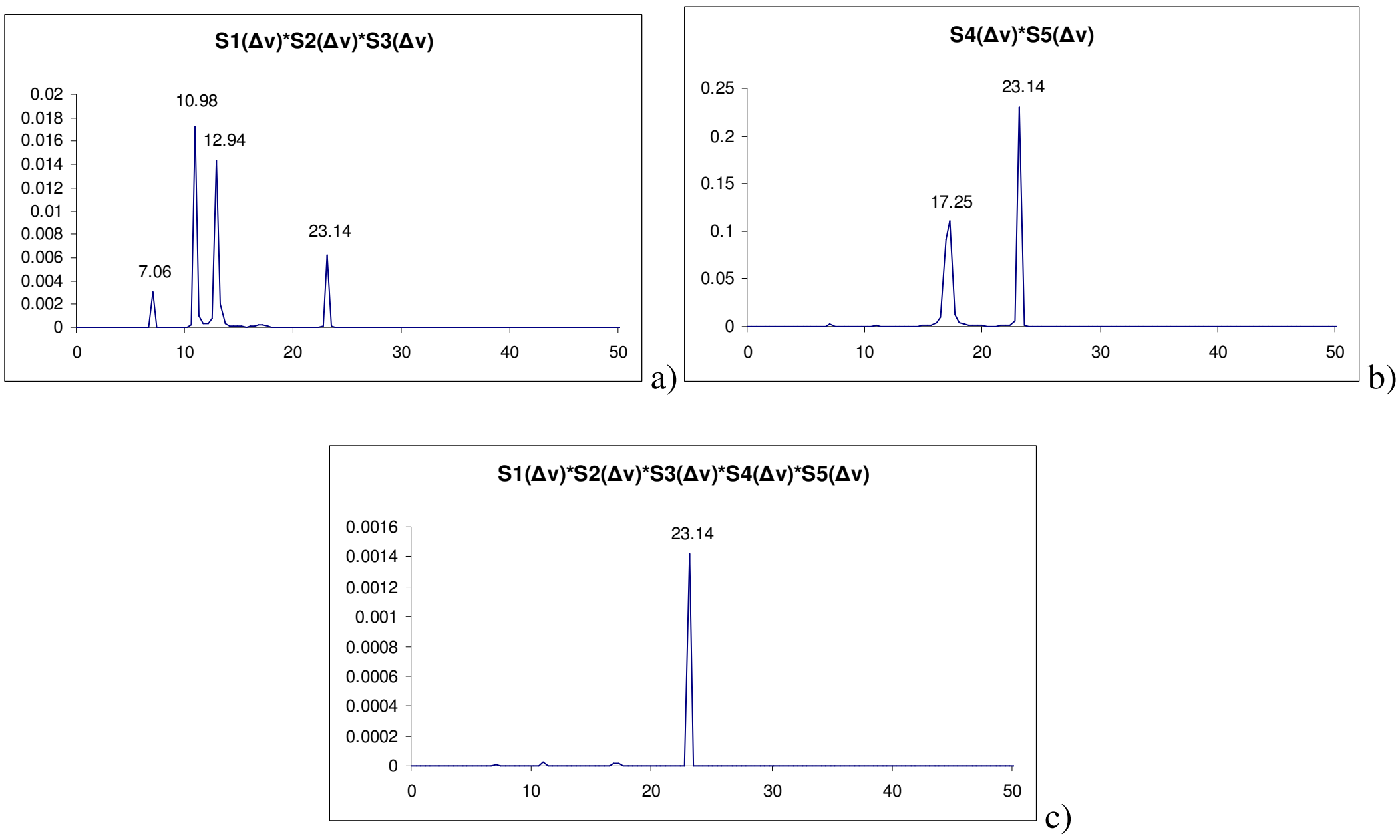

Figure 8. *The common frequencies*:
*a) emphasized in* $S_{123}(\Delta\nu) = S_1(\Delta\nu) \cdot S_2(\Delta\nu) \cdot S_3(\Delta\nu)$
*b) emphasized in* $S_{45}(\Delta\nu) = S_4(\Delta\nu) \cdot S_5(\Delta\nu)$
c) *emphasized in* $S_{12345}(\Delta\nu) = S_1(\Delta\nu) \cdot S_2(\Delta\nu) \cdot S_3(\Delta\nu) \cdot S_4(\Delta\nu) \cdot S_5(\Delta\nu)$

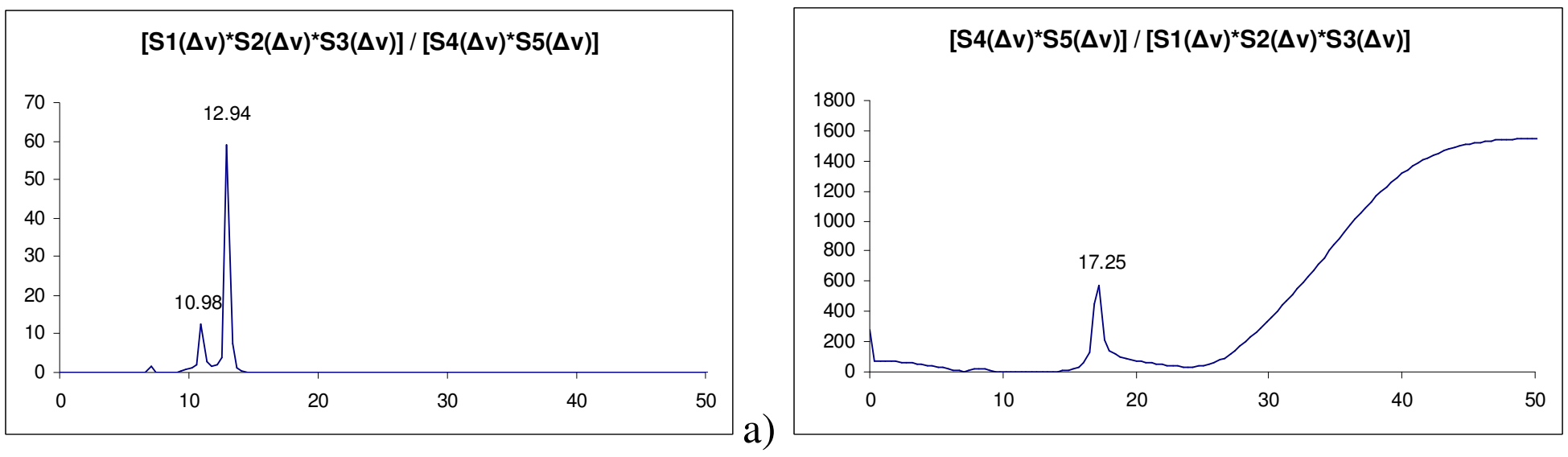

a) b)

Figure 9. *The non-common frequencies*:

*a) emphasized in* $S_{123}(\Delta\nu)\cdot S_{45}^{-1}(\Delta\nu)$ *b) emphasized in* $S_{45}(\Delta\nu)\cdot S_{123}^{-1}(\Delta\nu)$

The right-side "hill" from the Figure 9.b is the effect of the division between very small numbers such as $10^{-6}/10^{-9}$ , i.e. between spectral lines' amplitudes smaller than $10^{-6}$ obtained in the FFT analysis under the Microsoft Office Excel mathematical library, and can be eliminated by a simple stipulation on the arithmetic operation, as in Figure 10.

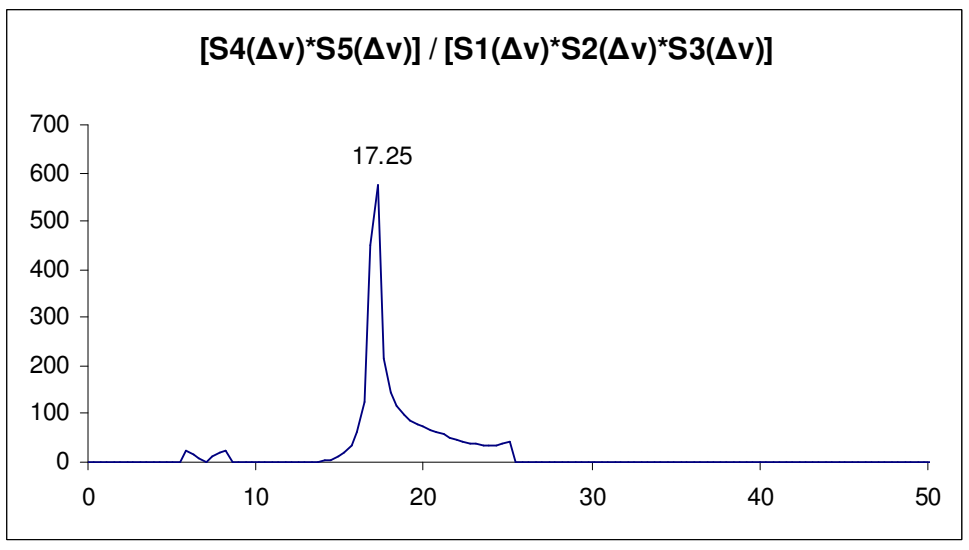

Figure 10. *The non-common frequencies*:

*emphasized in* $S_{45}(\Delta\nu)\cdot S_{123}^{-1}(\Delta\nu)$ *with conditioned division*

## 6. An experimental application

The above presented theoretical aspects were used in the analysis of the noise existing in the experimental measured signals of two absolute pressure transducers installed at the ends of a pipe through which water circulates, [11] – [12]. More precisely, in repeated experiments, for the two transducers' analyzed signals it was identified as significant common frequency a value $\nu_p$ very near-by the measured first modal frequency $\nu_1$ of the water-filled pipe's bending vibration (flow induced vibration), namely $|\nu_p - \nu_1|/\nu_1 < 1\%$ .

## 7. Conclusions

Starting with unconventional definitions for:

- *spectral line*: $(\nu, A(\nu))$;
- *(frequencies) discrete spectrum*: $S(\Delta\nu) = \{(\nu_i, A(\nu_i)) | i = 0,1,\cdots,N-1\}$, having the *resolution* $\Delta\nu = \nu_i - \nu_{i-1} = const.$, and $\nu_i = i\cdot\Delta\nu$ ;
- *congruent discrete spectra*: $S_1(\Delta\nu_1)$ and $S_2(\Delta\nu_2)$ if $\Delta\nu_1 = \Delta\nu_2 = \Delta\nu$ ;
- *correspondent spectral lines of two congruent discrete spectra*: $(\nu_k, A_1(\nu_k)) \in S_1(\Delta\nu)$ and $(\nu_k, A_2(\nu_k)) \in S_2(\Delta\nu)$;
- *the inverse of a (frequencies) discrete spectrum*: $S^{-1}(\Delta\nu) = \{(\nu_i, A^{-1}(\nu_i)) | i = 0,1,\ldots,N-1\}$ where $A^{-1}(\nu_i) = 1/A(\nu_i)$,

and considering the intuitive discernment criterion:

- *the more important frequency in signal has the greater magnitude among the spectral lines*,

the paper demonstrated the propositions:

*if* $S_1(\Delta\nu)$ *and* $S_2(\Delta\nu)$ *are two congruent discrete spectra, then:*

- *the discrete spectrum* $S_p(\Delta\nu) = S_1(\Delta\nu) \cdot S_2(\Delta\nu)$ *contains as emphasized the* ***common frequencies correspondent lines****, significantly existing in both* $S_1(\Delta\nu)$ *and* $S_2(\Delta\nu)$ *discrete spectra*;
- *the discrete spectrum* $S_r(\Delta\nu) = S_1(\Delta\nu) \cdot S_2^{-1}(\Delta\nu)$ *contains as emphasized the* ***non-common frequencie****s* ***correspondent lines****, existing significantly in* $S_1(\Delta\nu)$ *discrete spectrum, and insignificantly in* $S_2(\Delta\nu)$ *discrete spectrum,*

and, in general:

- *the discrete spectrum* $S_{\alpha\beta\cdots\zeta}(\Delta\nu) = S_\alpha(\Delta\nu) \cdot S_\beta(\Delta\nu) \cdot \ldots \cdot S_\zeta(\Delta\nu)$ *contains as emphasized the* ***common frequencies correspondent lines****, significantly existing in all congruent discrete spectra* $S_\alpha(\Delta\nu)$, $S_\beta(\Delta\nu)$,... $S_\zeta(\Delta\nu)$;
- *if* $S_{\alpha\beta\cdots\zeta}(\Delta\nu) = S_\alpha(\Delta\nu) \cdot S_\beta(\Delta\nu) \cdot \ldots \cdot S_\zeta(\Delta\nu)$ *and* $S_{ab\cdots z}(\Delta\nu) = S_a(\Delta\nu) \cdot S_b(\Delta\nu) \cdot \ldots \cdot S_z(\Delta\nu)$ *are congruent discrete spectra, then the discrete spectrum* $S_{\alpha\beta\cdots\zeta ab\cdots z}(\Delta\nu) = S_{\alpha\beta\cdots\zeta}(\Delta\nu) \cdot S_{ab\cdots z}^{-1}(\Delta\nu)$ *contains as empha-sized the* ***non-common frequencies correspondent lines****, existing significantly in all* $S_\alpha(\Delta\nu)$, $S_\beta(\Delta\nu)$,... $S_\zeta(\Delta\nu)$ *discrete spectra, and insignificantly in all* $S_a(\Delta\nu)$, $S_b(\Delta\nu)$,... $S_z(\Delta\nu)$ *discrete spectra,* statements which allow the development of a practical method dedicated to detection of the *common* and *non-common frequencies* from two or more *congruent discrete spectra*.

Mathematically funded and fully illustrated by numerical simulations made using the Microsoft Office Excel computational tools, the proposed method can be implemented, at an easy rate, in any automated signals monitoring task (i.e.. using a PC applications), with an adequate degree of reliability.

## R E F E R E N C E S

[1] ***, *ELECTRICAL MEASUREMENT, SIGNAL PROCESSING, and DISPLAYS*, Edited by JOHN G. WEBSTER, Copyright © 2004 by CRC Press LLC

[2] Bernd Girod, *Signals and Systems*, Copyright © 2001 by John Wiley & Sons Ltd.

[3] E.S. Gopi, *Algorithm Collections for Digital Signal Processing Applications Using Matlab*®, All Rights Reserved © 2007 Springer

[4] Robert M. Gray and Lee D. Davisson, *An Introduction to Statistical Signal Processing*, Copyright © 2004 by Cambridge University Press.

[5] Steven T. Karris, *Signals and Systems with MATLAB® Computing and Simulink® Modeling*, Copyright © 2008 Orchard Publications

[6] D.S.G. Pollock, *A Handbook of Time-Series Analysis, Signal Processing and Dynamics*, Copyright © 1999 by ACADEMIC PRESS

[7] Peter Seibt, *Algorithmic Information Theory, Mathematics of Digital Information Processing*, Copyright © Springer-Verlag Berlin Heidelberg 2006

[8] Steven W. Smith, *The Scientist and Engineer's Guide to Digital Signal Processing*, Second Edition, Copyright © 1997-1999 by California Technical Publishing

[9] R.W. Stewart and M.W. Hoffman, *Digital Signal Processing, An "A" to "Z"*, Copyright © BlueBox Multimedia, R.W. Stewart 1998

[10] Li Tan, *Digital Signal Processing, Fundamentals and Applications*, Copyright © 2008, Elsevier Inc.

[11] C. Doca. C. Păunoiu, I. Tebeică, N. Anghel, I. Costache and C. Tudor, *MONITORING THE PIPE'S BENDING VIBRATIONS USING ABSOLUTE PRESSURE TRANSDUCERS?*, The International Symposium for Nuclear Energy SIEN 2013 "Nuclear Power – The Today's Challenge", 10-14 November, 2013, Bucharest, Romania

[12] C. Doca and C. Păunoiu, *AN INTUITIVE METHOD TO AUTOMATICALLY DETECT THE COMMON AND NON-COMMON FREQUENCIES FOR TWO OR MORE TIME-VARYING SIGNALS,* Scientific Bulletin of "Politehnica" University of Bucharest, Series C, Vol. 76, Iss. 4, 2014, ISNN 2286-3540, pag. 267-278